\documentclass[
  aps,
  prx,
  twocolumn,
  english,
  superscriptaddress,
  floatfix,
  longbibliography,
  10pt
]{revtex4-2}

\usepackage{amsmath,amssymb,amsthm}
\usepackage{latexsym}
\usepackage{indentfirst}
\usepackage{graphicx}

\usepackage{placeins}
\usepackage{booktabs}
\usepackage{algorithm}
\usepackage{algorithmic}
\usepackage{multirow}
\usepackage[caption=false]{subfig}
\usepackage{graphicx,color}
\usepackage{diagbox}
\usepackage{braket}
\usepackage[normalem]{ulem}
\usepackage{hyperref}
\usepackage{cleveref}
\usepackage{newfloat}

\theoremstyle{plain}
\newtheorem{theorem}{Theorem}
\newtheorem{lemma}[theorem]{Lemma}

\newtheorem{prop}[theorem]{Proposition}

\theoremstyle{definition}

\theoremstyle{remark}

\newtheoremstyle{cited}%
  {3pt}%
  {3pt}%
  {\itshape}%
  {}%
  {\bfseries}%
  {.}%
  {.5em}%
  {\thmname{#1} \thmnumber{#2} \thmnote{\normalfont#3}}%

\theoremstyle{cited}

\usepackage{amsmath,amsfonts,bm}

\def\eqref#1{equation~\ref{#1}}

\def\1{\bm{1}}

\DeclareMathAlphabet{\mathsfit}{\encodingdefault}{\sfdefault}{m}{sl}
\SetMathAlphabet{\mathsfit}{bold}{\encodingdefault}{\sfdefault}{bx}{n}

\newcommand{\R}{\mathbb{R}}
\newcommand{\C}{\mathbb{C}}

\newcommand{\KL}{D_{\mathrm{KL}}}

\DeclareMathOperator*{\argmin}{arg\,min}

\newcommand{\bbm}{\begin{bmatrix}}
\newcommand{\ebm}{\end{bmatrix}}
\newcommand{\DMRG}{\mathrm{cross}}
\newcommand{\pcausal}{p^{\mathrm{causal}}}
\newcommand{\gradientat}[2]{\left.#1\right|_{#2}}

\renewcommand{\KL}[2]{D_{\mathrm{KL}}\left(#1 \parallel #2\right)}

\newcommand{\NLL}{\mathcal{L}}
\newcommand{\PTS}{\hat p_{\mathrm{TS}}}

\definecolor{fern}{RGB}{79,143,0}

\usepackage{mathtools}
\usepackage{etoolbox}
\DeclareFontFamily{U}{matha}{\hyphenchar\font45}
\DeclareFontShape{U}{matha}{m}{n}{
      <5> <6> <7> <8> <9> <10> gen * matha
      <10.95> matha10 <12> <14.4> <17.28> <20.74> <24.88> matha12
      }{}
\DeclareSymbolFont{matha}{U}{matha}{m}{n}
\DeclareFontSubstitution{U}{matha}{m}{n}

\DeclareFontFamily{U}{mathx}{\hyphenchar\font45}
\DeclareFontShape{U}{mathx}{m}{n}{
      <5> <6> <7> <8> <9> <10>
      <10.95> <12> <14.4> <17.28> <20.74> <24.88>
      mathx10
      }{}
\DeclareSymbolFont{mathx}{U}{mathx}{m}{n}
\DeclareFontSubstitution{U}{mathx}{m}{n}

\DeclareMathDelimiter{\vvvert}{0}{matha}{"7E}{mathx}{"17}
\DeclarePairedDelimiterX{\normi}[1]
  {\vvvert}
  {\vvvert}
  {\ifblank{#1}{\:\cdot\:}{#1}}

\begin{document}

\title{Initialization and training of matrix product state probabilistic models}
\author{Xun Tang}
\affiliation{Department of Mathematics, Stanford University}
\author{Yuehaw Khoo}
\affiliation{Department of Statistics, University of Chicago}
\author{Lexing Ying}
\affiliation{Department of Mathematics, Stanford University}
\date{\today}

\begin{abstract}
    Modeling probability distributions via the wave function of a quantum state is central to quantum-inspired generative modeling and quantum state tomography (QST). We investigate a common failure mode in training randomly initialized matrix product states (MPS) using gradient descent. The results show that the trained MPS models do not accurately predict the strong interactions between boundary sites in periodic spin chain models. In the case of the Born machine algorithm, we further identify a \emph{causality trap}, where the trained MPS models resemble causal models that ignore the non-local correlations in the true distribution.
    We propose two complementary strategies to overcome the training failure---one through optimization and one through initialization.
    First, we develop a natural gradient descent (NGD) method, which approximately simulates the gradient flow on tensor manifolds and significantly enhances training efficiency. Numerical experiments show that NGD avoids local minima in both Born machines and in general MPS tomography. Remarkably, we show that NGD with line search can converge to the global minimum in only a few iterations.
    Second, for the BM algorithm, we introduce a warm-start initialization based on the TTNS-Sketch algorithm. We show that gradient descent under a warm initialization does not encounter the causality trap and admits rapid convergence to the ground truth.
\end{abstract}
\maketitle

\section{Introduction}

Recent advances in generative modeling enable one to learn complex high-dimensional distributions
\cite{larochelle2011neural,germain2015made,salimans2017pixelcnn++,kingma2013auto,lecun2006tutorial, rezende2015variational,goodfellow2014generative,song2020score}.
For discrete distributions, tensor network (TN) architectures have emerged as a prevalent method for probabilistic modeling \cite{han2018unsupervised, hur2022generative,tang2023generative,novikov2021tensor}. One prominent class of tensor network models takes the perspective of Born machine (BM) \cite{han2018unsupervised}. Under this setting, one takes a tensor network to represent a quantum state. The probabilistic function follows from the Born rule and is thus represented by the squared modulus of the quantum state. 

This work focuses on Born machines using matrix product states (MPS)---a one-dimensional tensor network with strong expressivity \cite{white1992density}.
Qubits entangled with a specific arrangement of local quantum circuits can be exactly modeled by an MPS representation \cite{glasser2019expressive}. The BM algorithm minimizes the negative log-likelihood (NLL) over observed samples, fitting the squared MPS amplitude to empirical data. The BM formulation is a special case of MPS quantum state tomography \cite{cramer2010efficient,lanyon2017efficient}, with the key distinction that the BM algorithm fits MPS models against measurements made on the computational basis. Training of the BM algorithm can be done classically, involving only conventional numerical linear algebra. In contrast, explicitly training a variational quantum circuit to fit the given samples would involve automatic differentiation \cite{mitarai2018quantum, bergholm2018pennylane, schuld2019evaluating} on quantum hardware, which is more difficult due to the presence of noise in obtaining the gradient through measurement.

The use of MPS is proven successful in minimizing the energy of a quantum system, a well-celebrated example being the density-matrix renormalization group (DMRG) algorithm \cite{white1992density}. However, unsupervised generative modeling with MPS is a nonconvex optimization setting with unique challenges. While there is extensive literature on training neural networks in such settings (e.g., \cite{glorot2010understanding,pascanu2013difficulty,erhan2010does}), the performance of the MPS ansatz in general optimization tasks is less understood.

In this work, we use numerical evidence to show that gradient descent (GD) can lead to training failures for nonconvex optimization tasks under the MPS ansatz, both for the BM algorithm and for MPS tomography in general. 
We demonstrate that a randomly initialized MPS model fails to converge to the global minimum using standard gradient descent. 
We show that even substantial over-parameterization does not enable the model to escape these minima. The trained model exhibits rank degeneracy and overlooks important non-local correlations, making inferences from such models questionable. In the BM setting, we characterize the failure mode as a ``causality trap," wherein training converges to a simplified causal model, a graphical model class with significantly less approximation power. In addition, the causality trap phenomenon also occurs in the DMRG-type training method considered in \cite{han2018unsupervised}.

The observed local minima issues share some similarities with the barren plateau phenomenon in quantum machine learning \cite{mcclean2018barren,marrero2021entanglement,cerezo2021cost,
cerezo2021higher,anschuetz2022beyond,anschuetz2021critical,dborin2022matrix}. However, we show the challenge facing an MPS model in this case is a mild local minimum problem typical for nonconvex optimization. This work proposes two training strategies to prevent such issues.

First, we propose a natural gradient descent (NGD) method that performs optimization in the function space of high-dimensional tensors, rather than directly on parameters. Mathematically, the proposed training process is the discretization of a projected gradient flow in the space of quantum states. The proposal allows efficient convergence to the global minimum in both BM and MPS tomography.

Second, for the BM algorithm setting, we propose a warm-start initialization protocol. We propose a warm-start procedure using the tensor tree network states via sketching (TTNS-Sketch) algorithm, which gives a consistent density estimator with sample complexity guarantees \cite{tang2023generative}. Our finding shows that a non-random initialization allows the gradient descent method to converge to the ground truth.

\begin{figure}
    \centering
    \includegraphics[width=0.9\columnwidth]{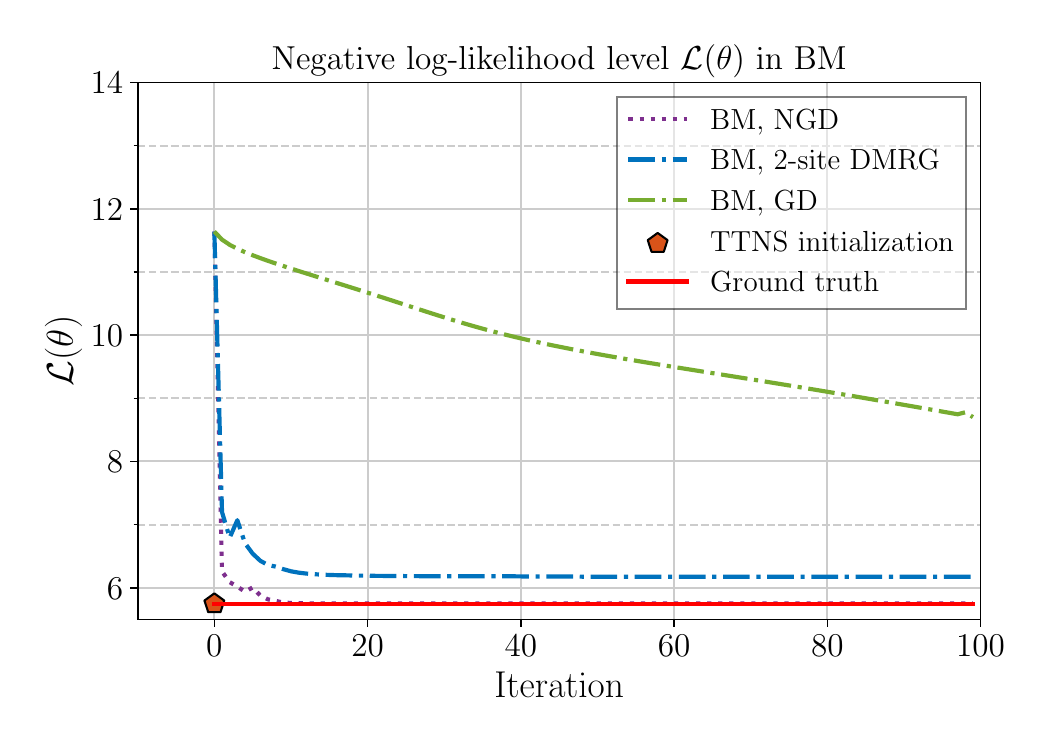}
    \caption{Illustration of NGD and warm initialization for the Born machine algorithm. One sees that the gradient descent method and the 2-site DMRG method do not converge to the optimal log-likelihood level. The ground truth model is a periodic ferromagnetic Ising model, and the experiment details are in \Cref{sec: analysis}.}
    \label{fig:combined_NLL_plot_BM}
\end{figure}

The effect of adopting the proposed methods is shown in \Cref{fig:combined_NLL_plot_BM}. 
For the BM algorithm, the choices of initialization are (i) random initialization and (ii) the proposed warm initialization based on TTNS-sketch. The choices of training methods include: (a) the GD method, (b) the 2-site DMRG method in \cite{han2018unsupervised}, and (c) the proposed NGD method. 
In terms of training method, \Cref{fig:combined_NLL_plot_BM} shows that the NGD method can converge to the global minimum under random initialization, while both GD and 2-site DMRG encounter local minima issues. Moreover, in \Cref{fig:combined_NLL_plot_BM}, the warm initialization is sufficiently close to the global minimum that no further training is needed. Therefore, one can see that the local minima issue in the Born machine can be addressed through either improved initialization or improved training methods.

The remainder of the paper is organized as follows: \Cref{sec: analysis} presents the causality trap in Born machine settings; \Cref{sec: local minima in MPS tomography} reports local minima issues in MPS tomography; \Cref{sec: ngd} presents the NGD method for MPS optimization; \Cref{sec: TTNS initialization} introduces the TTNS-Sketch warm initialization for the BM algorithm; and \Cref{sec: conclusion} offers concluding remarks.

\begin{figure}
    \centering
    \subfloat[Cycle graph with \(n = 16\)]{
        \includegraphics[width=0.4\columnwidth]{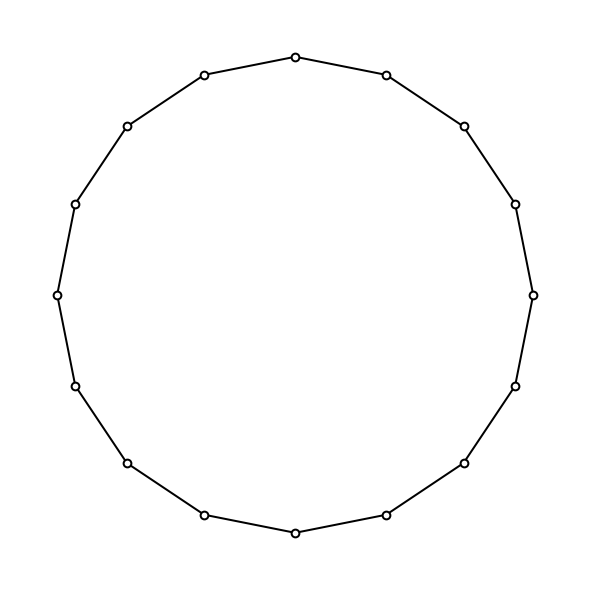}
        \label{fig:circle}
    }
    \hspace{0.05\columnwidth}
    \subfloat[Line graph with \(n = 16\)]{
        \includegraphics[width=0.4\columnwidth]{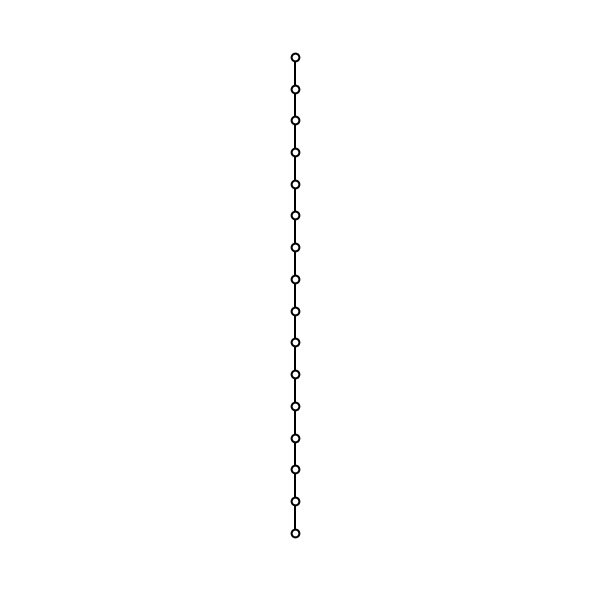}
        \label{fig:line_16}
    }
    \caption{Graphical representations of the underlying model \(p^{\star}\) (Fig.~\ref{fig:circle}) and of the mis-specified model \(\pcausal\) (Fig.~\ref{fig:line_16}).}
    \label{fig:graph plot}
\end{figure}

\section{Causality trap in the Born machine algorithm}\label{sec: analysis}

\subsection{Background in Born machine}
We begin with a brief introduction to maximum likelihood estimation and the Born machine (BM) algorithm. Suppose one is given an \emph{underlying distribution} \(p^{\star}\) and a parameterized family of distributions \(\{p_{\theta}\}_{\theta \in \Theta}\). Given a dataset \(\mathcal{T}\) of samples drawn from \(p^{\star}\), the negative log-likelihood (NLL) measures how a parameterized density fits \(\mathcal{T}\), and is defined as follows:
\begin{equation}\label{eqn: NLL formula}
    \NLL_{\mathrm{BM}}(\theta) = -\frac{1}{|\mathcal{T}|}\sum_{y \in \mathcal{T}}{\log(p_{\theta}(y))},
\end{equation}
and we write \(\NLL = \NLL_{\mathrm{BM}}\) for the remainder of this section.

The Born machine uses a matrix product state as the tensor network ansatz. In particular, the parameter \(\theta = (G_{k})_{k = 1}^{n}\) is a collection of tensor components, where \(G_{1} \in \R^{ 2 \times r_{1}}\), \(G_{i} \in \R^{r_{i-1} \times 2 \times r_{i}}\) for \(i = 2,\ldots, n-1\), and \(G_{n} \in \R^{r_{n-1} \times 2 }\). The MPS \(q_{\theta}\) takes the following form:
\begin{equation}
\begin{aligned}
    &q_{\theta}(x_{1}, \ldots, x_{n}) \\ =&\sum_{\alpha_{1}, \ldots, \alpha_{n-1}}G_{1}(x_{1}, \alpha_{1})G_{2}(\alpha_{1}, x_{2}, \alpha_{2})\cdots G_{n}(\alpha_{n-1}, x_{n}),
\end{aligned}
\end{equation}
and the associated equation for the density function \(p_{\theta}\) is
\begin{equation}\label{eq: BM equation}
p_{\theta}(x) =\frac{|q_{\theta}(x)|^2}{Z_{\theta}},
\end{equation}
where \(Z_{\theta} = \sum_{z} |q_{\theta}(z)|^2\) is the associated normalization constant and can be efficiently computed by applying tensor contractions. The goal is to find \(\hat{\theta} = \argmin_{\theta} \NLL(\theta)\), and the resultant \(p_{\hat{\theta}}\) is the maximum likelihood estimator of the density \(p^{\star}\).

\begin{figure*}
    \centering
         \includegraphics[width=1.8\columnwidth]{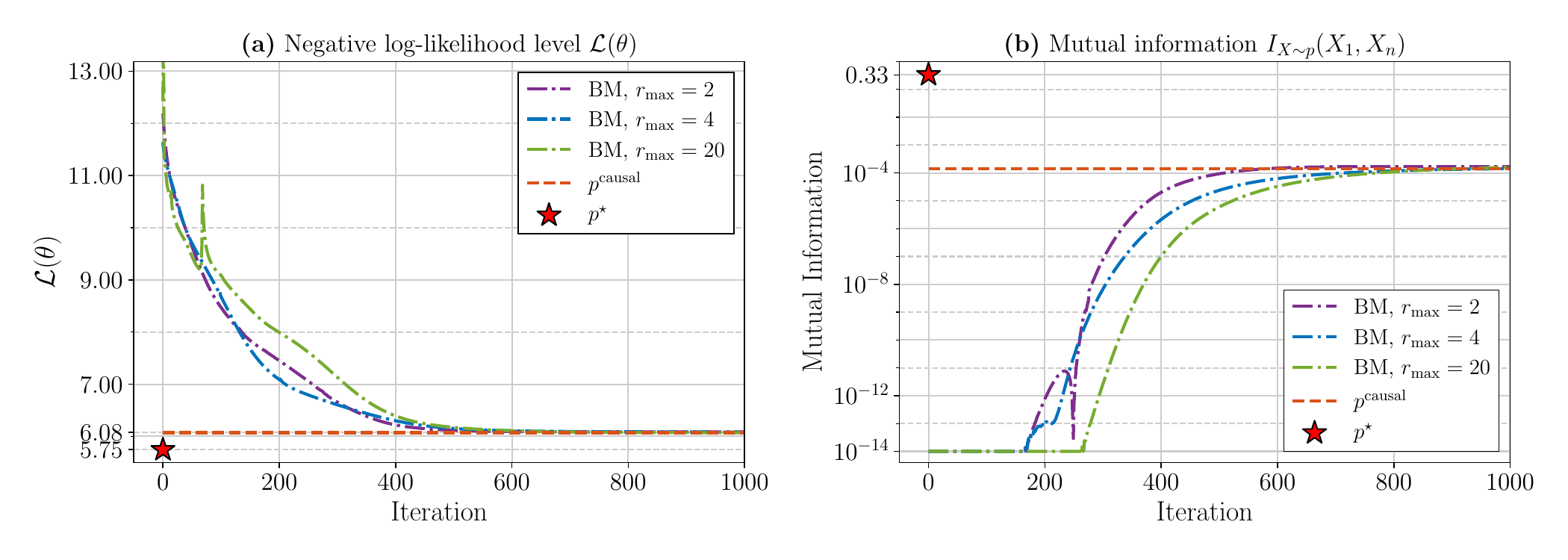}
    \caption{Performance of Born machine algorithm for the periodic spin system model in \Cref{eqn: periodic model}. The models are initialized randomly and trained under gradient descent. The NLL gap is \(0.33\), which coincides with the mutual information level of \((X_1, X_n)\) in \(p^{\star}\). \Cref{appendix: KL analysis} shows that the NLL gap and the mutual information level are approximately equal under the causality trap.}
    \label{fig:result of GD}
\end{figure*}

\begin{figure*}
    \centering
    \includegraphics[width=1.8\columnwidth]{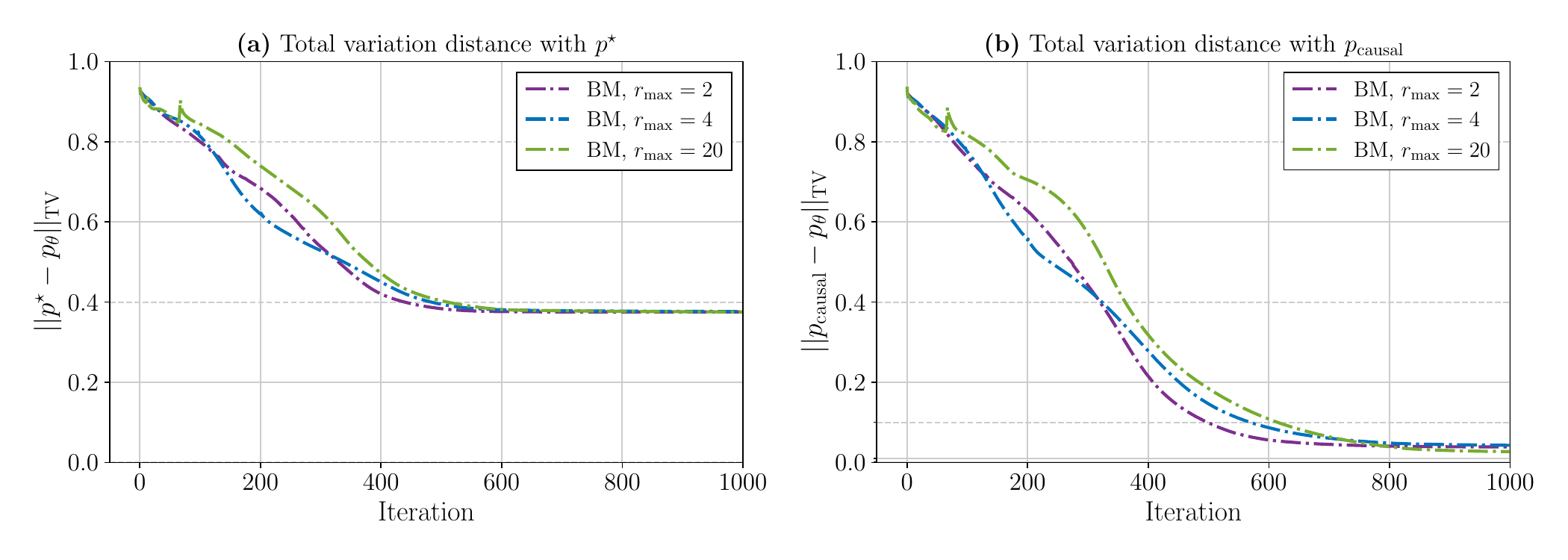}
    \caption{Plot of total variation (TV) distance of the trained Born machine model with respect to the true model \(p^{\star}\) and the causal model \(\pcausal\). The setting is the same as in \Cref{fig:result of GD}. One can see that the trained BM model is much closer to the causal model than to the true model. The TV distance is defined by \(\lVert p -p' \rVert_{\mathrm{TV}} = \frac{1}{2}\lVert p - p' \rVert_{1}\).}
    \label{fig:TV of GD}
\end{figure*}

\subsection{Causality trap under a periodic Ising model}\label{sec: experiment setting}

We consider a simple distribution given by the ferromagnetic Ising model with a periodic boundary condition:
\begin{equation}\label{eqn: periodic model}
    p^{\star}(x_{1},\ldots, x_{n}) \propto \exp{\left(-\beta\sum_{(i,j) \in \mathrm{cycle}(n)}x_{i} x_{j}\right)},
\end{equation}
where \(x_{i} \in \{-1,1\}\) and \(\mathrm{cycle}(n)\) is the cycle graph over \(n\) sites.
The model in \Cref{eqn: periodic model} can be characterized by a graphical model over a cycle; see \Cref{fig:circle}. For a measure for model complexity, we define \emph{maximal internal rank} \(r_{\text{max}} := \max_{i}r_{i}\), where \(\{r_i\}_{i = 1}^{n-1}\) is the internal rank determining the size of each tensor component \(G_i\). The model in \Cref{eqn: periodic model} can be represented by a BM ansatz with \(r_{\text{max}} \leq 4\).

The training process is done by gradient descent, and the results are obtained using the existing algorithmic implementation from \cite{glasser2019expressive}.
For the experiment, we let \(\beta = 1\), \(n = 16\). We select a large sample size by taking \(\lvert \mathcal{T} \rvert = 2^{15}\) for training.  
The training result is plotted in \Cref{fig:result of GD}.
For all choices of parameter sizes, the learned model stays at a sub-optimal NLL level with a significant gap from the global minimum. One sees that the NLL gap persists even under the over-parameterization setting of \(r_{\text{max}} = 20\).

Moreover, the learned BM model fails to model the important boundary correlation. For \(X = (X_1, \ldots, X_n) \sim p^{\star}\), the mutual information for the variable pair \((X_1, X_n)\) is large in \(p^{\star}\). However, as shown in \Cref{fig:result of GD}, the trained BM models predict a weak mutual information level, and so the trained model fails to capture the interaction between the variable pair \((X_1, X_n)\).

The \emph{causality trap} refers to the phenomenon that the training dynamics of BM effectively converge to the following causal model
\begin{equation}\label{eqn: best causal model}
    \pcausal(x_{1},\ldots, x_{n}) \propto \exp{\left(-\sum_{(i,j) \in \mathrm{path}(n)}x_{i} x_{j}\right)},
\end{equation}
where \(\mathrm{path}(n)\) is the path graph over \(n\) sites. The model in \Cref{eqn: best causal model} is similar to the one in \Cref{eqn: periodic model} but misses the important interaction term coming from the edge linking site \(1\) and site \(n\). As can be seen in \Cref{fig:result of GD}, the trained model under gradient descent closely matches \(\pcausal\) in both NLL and in the mutual information for \((X_1, X_n)\). Furthermore, we show in \Cref{fig:TV of GD} that the trained BM model is very close to \(\pcausal\) in terms of the total variation (TV) distance.

We remark that \(\pcausal\) is representable by a BM ansatz requiring only \(r_{\text{max}} = 2\). Therefore, under gradient descent, the training dynamics of BM favors outputting rank degenerate models, even though the internal rank is typically set large to ensure approximation power. As seen in \Cref{fig:result of GD}, one can see that the causality trap occurs even when \(r_{\text{max}} = 20\), which shows that the causality trap persists even under substantial over-parameterization. 

For larger \(n\), evaluating the causality trap through the TV distance has an \(O(2^n)\) scaling. In \Cref{appendix: KL analysis}, we perform a detailed analysis of the causality trap and show that the causality trap can be exactly characterized by the trained model matching \(\pcausal\) in NLL and the mutual information for \((X_1, X_n)\).

\section{Local minima in MPS state tomography}\label{sec: local minima in MPS tomography}
\subsection{Background in MPS tomography}
Quantum state tomography (QST) is the task of finding a quantum state from measurement outcomes \cite{vogel1989determination, leonhardt1995quantum, white1999nonmaximally, cramer2010efficient, torlai2018neural}. We take the \(n\)-bit setting in this section for simplicity. One has \(B\) copies of a ground truth quantum state \(\ket{\phi}\). In this case, one is given \(B\) unitary transformations \(\{U^{(j)} \in U(2^n)\}_{j = 1}^{B}\). For \(j = 1, \ldots, B\), one performs a computational-basis measurement on \(U^{(j)}\ket{\phi}\) and receives a measurement outcome \(\ket{b^{(j)}} \in \{0, 1\}^{n}\). The input dataset to the learning task is \(\mathcal{T} = \{(b^{(j)}, U^{(j)})\}_{j = 1}^{B}\). 
QST is typically done by maximum likelihood. For a parameterized quantum state \(\ket{\psi_{\theta}}\), one minimizes the NLL defined as follows
\[
    \NLL_{\mathrm{QST}}(\theta) = -\frac{1}{|\mathcal{T}|}\sum_{b, U \in \mathcal{T}}{\log\left(\lvert\bra{b} U\ket{\psi_{\theta}}\rvert^2\right)}.
\]
Similar to the BM case, the goal is to find \(\hat{\theta} = \argmin_{\theta}\NLL_{\mathrm{QST}}(\theta)\), and \(\ket{\psi_{\hat{\theta}}}\) is the maximum likelihood estimator one wishes to obtain.

The MPS tomography assumes a complex MPS ansatz to model \(\ket{\psi}\). In this case, one uses the parameter \(\theta\) to encode tensor components \((G_{k})_{k = 1}^{n}\), where \(G_{1} \in \C^{ 2 \times r_{1}}\), \(G_{i} \in \C^{r_{i-1} \times 2 \times r_{i}}\) for \(i = 2,\ldots, n-1\), and \(G_{n} \in \C^{r_{n-1} \times 2 }\). The MPS state \(\ket{\psi_{\theta}} \in \C^{2^n}\) satisfies the following equation:
\begin{equation*}
\begin{aligned}
    &\ket{\psi_{\theta}}_{x_1, \ldots, x_n} \\ =&\frac{1}{Z}\sum_{\alpha_{1}, \ldots, \alpha_{n-1}}G_{1}(x_{1}, \alpha_{1})G_{2}(\alpha_{1}, x_{2}, \alpha_{2})\cdots G_{n}(\alpha_{n-1}, x_{n}),
\end{aligned}
\end{equation*}
where \(Z_{\theta} = \braket{\psi_{\theta} |\psi_{\theta}}\) is the normalization constant. 

In particular, the BM algorithm is equivalent to a special case of MPS tomography where the state \(\ket{\phi}\) undergoes computational-basis measurement without applying unitary transformations.

\begin{figure}
    \centering
         \includegraphics[width=0.9\columnwidth]{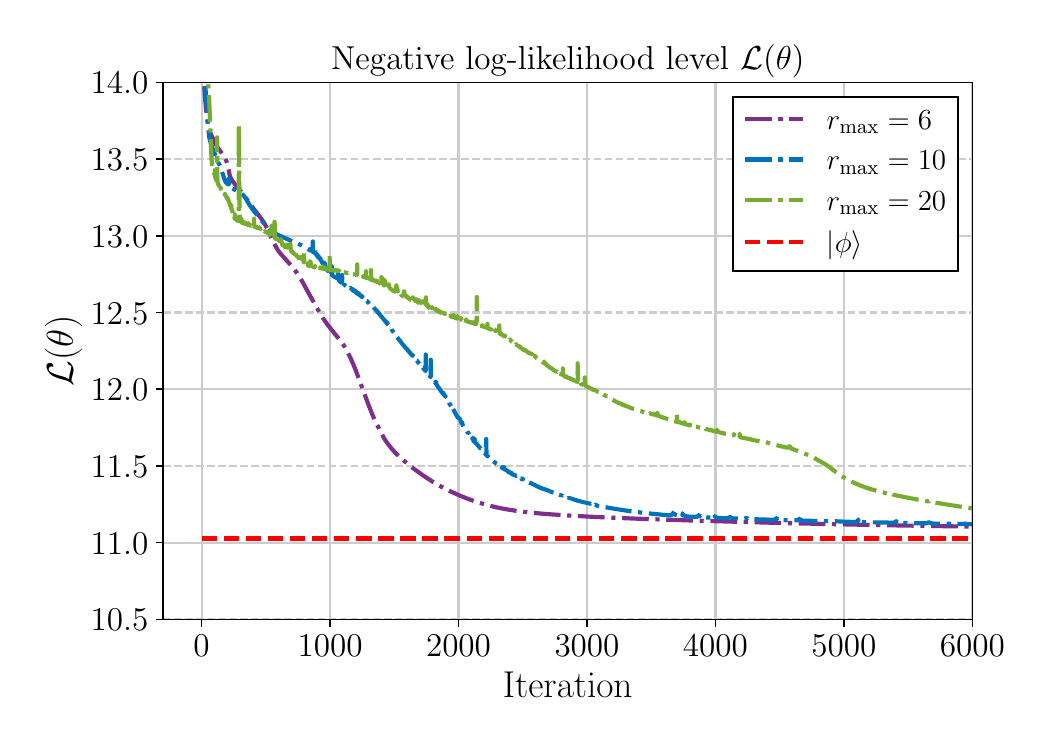}
    \caption{Performance of MPS tomography algorithm for the ground state of the periodic TFIM model in \Cref{eqn: TFIM}. The models are initialized randomly and trained under gradient descent.}
    \label{fig:result of GD in MPS}
\end{figure}
\subsection{Local minima under a periodic TFIM model}\label{sec: experiment setting TFIM}
We consider the task of quantum state tomography for the ground state of the 1D ferromagnetic transverse field Ising model (TFIM). The Hamiltonian of the TFIM model is
\begin{equation}\label{eqn: TFIM}
    H = -J \sum_{(i,j) \in \mathrm{cycle}(n)} \sigma_{i}^{Z} \sigma_{j}^{Z} -  h \sum_{i} \sigma_{i}^{X},
\end{equation}
where \(\sigma_i^{Z}\) (resp. \(\sigma_i^{X}\)) is the Pauli-Z (resp. Pauli-X) matrix on site \(i\). We consider a critical point by taking \(J = h = 1\). We obtain the ground state \(\ket{\phi}\) as an MPS by using the density matrix renormalization group (DMRG) algorithm. In particular, we use DMRG to model the ground state as an MPS \(\ket{\phi}\) of maximal internal bond dimension \(r_{\mathrm{max}} = 6\). 

For the experiment, we take a system size of \(n = 20\). We record \(B = 20000\) samples of \(\ket{\phi}\) by random Pauli measurements. In other words, for each \(j = 1, \ldots, B\) and \(i = 1, \ldots, n\), we select \(U^{(j)}_{i}\) to be a unitary matrix on site \(i\), and we uniformly choose between the \(X, Y, Z\) basis on each site. Then, we take each unitary transformation \(U^{(j)}\) to be \(U^{(j)} = \bigotimes_{i=1}^{n}U^{(j)}_{i}\). We remark that the state \(U^{(j)}\ket{\phi}\) is also an MPS of the same shape as \(\ket{\phi}\), and so performing the computational-basis measurement can be done classically.

Similarly to the BM case, we use the gradient descent algorithm to perform training, and the tensor components of the MPS models are randomly initialized. The training result is in Figure \ref{fig:result of GD in MPS}. For all choices of parameter sizes, the learned model stays at a sub-optimal NLL level with a significant gap from the NLL of the true state \(\ket{\phi}\). One sees that the NLL gap persists even under the over-parameterization setting of \(r_{\text{max}} = 20\).

\section{Natural gradient descent algorithm for MPS optimization}\label{sec: ngd}
In this section, we propose a natural gradient descent (NGD) method, which improves on the gradient descent (GD) method used for the MPS ansatz.

We explain the main idea of the NGD method for MPS. In the general variational setting, one has a parameterized MPS family \(q_{\theta}\) with tensor components \(\theta = (G_{k})_{k=1}^{n}\). The goal is to minimize a loss function \(\mathcal{L}(\theta)\) defined on the parameter space. The NGD method can be summarized as the following optimization task:
\begin{equation}\label{eqn: NGD for TN}
\begin{aligned}
    \theta_{t+1} = \theta_{t} +  \argmin_{\delta\theta} \left<\nabla_{\theta}\mathcal{L}|_{\theta = \theta_{t}}, \delta\theta \right> +  \frac{1}{2}\eta \lVert q_{\theta_{t} + \delta \theta} - q_{\theta_{t}} \rVert_{F}^2,
\end{aligned}
\end{equation}
where \(\lVert \cdot \rVert_{F}\) denotes the Frobenius norm.

We discuss the difference between GD and NGD. In \Cref{eqn: NGD for TN}, if one replaces \(\frac{1}{2}\eta \lVert q_{\theta_{t} + \delta \theta} - q_{\theta_{t}} \rVert_{F}^2\) with \(\frac{1}{2}\eta \lVert \delta \theta \rVert_{F}^2\), then one would recover the GD algorithm. The learning rate is \(\eta^{-1}\) for both cases. Essentially, both algorithms consider the minimization of the linear approximation of \(\mathcal{L}(\theta)\) around \(\theta = \theta_t\), but NGD uses \(\frac{1}{2}\eta \lVert q_{\theta_{t} + \delta \theta} - q_{\theta_{t}} \rVert_{F}^2\) as the curvature term for regularization. The main benefit of the NGD approach is that its curvature term considers the variation in the exponential-sized tensor space instead of the parameter space. In Appendix \ref{sec: multi-linearity induced vanishing gradient}, we give a toy example in which NGD ensures training success, whereas GD experiences a vanishing/exploding gradient problem. One additional benefit is that the NGD approach is independent of the gauge degree of freedom of the MPS. 

\subsection{Algorithm}
In practice, the NGD method is implemented with a sequential tensor component update. Writing \(\delta \theta = (\delta G_1, \ldots, \delta G_{n})\), one can check that the minimization task in \Cref{eqn: NGD for TN} is a quadratic program in each \(\delta G_i\) for \(i \in \{1, \ldots, n\}\). Therefore, to update \(\theta_{t}\), one picks a site \(i\) and performs the optimization task over \(\delta G_i\) in \Cref{eqn: NGD for TN}, and the update can be done analytically.

Our proposed NGD procedure is summarized in \Cref{alg:1}. Due to the 1D geometry of MPS, the sequence of site-wise update is most efficient if one performs a forward sweep (picking \(i\) from \(1\) to \(n\)) followed by a backward sweep (picking \(i\) from \(n\) to \(1\)). The reason for the site update schedule is to cache and reuse intermediate tensor components for optimal efficiency. When \(\mathcal{L}\) is the NLL loss, for example, running \Cref{alg:1} has only a time complexity of \(O(n)\).

\begin{algorithm}[H]
\caption{Natural gradient descent update with optional line search.}
\label{alg:1}
\begin{algorithmic}[1]
\REQUIRE Loss function \(\mathcal{L}\).
\REQUIRE Current tensor component \(\theta_{t}\), parameter \(\eta\).

\FOR{\(i = 1,\ldots, n\) \textbf{and then} \(i = n,\ldots,1\)}
\STATE \(S_{i} \gets \{(\delta G_{k})_{k = 1}^{n} \mid \delta G_{k} = 0 \, \forall \, k \not = i\}\)
\STATE 
\(
\delta \theta_{t} \gets \argmin_{\delta\theta \in S_{i}} \left<\nabla_{\theta}\mathcal{L}|_{\theta = \theta_{t}}, \delta\theta \right> + \frac{1}{2}\eta \lVert q_{\theta_{t} + \delta \theta} - q_{\theta_{t}} \rVert_{F}^2
\)
\IF{using line search}
    \STATE Find \(\alpha = \argmin_{\alpha > 0} \mathcal{L}(\theta_t + \alpha \delta \theta_{t})\)
    \STATE Update \(\theta_{t} \gets \theta_{t} + \alpha \delta \theta_{t}\)
\ELSE
    \STATE Update \(\theta_{t} \gets \theta_{t} + \delta \theta_{t}\)
\ENDIF
\ENDFOR

\STATE Set \(\theta_{t+1} \gets \theta_{t}\).

\RETURN \(\theta_{t+1}\)
\end{algorithmic}
\end{algorithm}

In Appendix \ref{appendix: NGD versus DMRG}, we show that the NGD step in \Cref{alg:1} can be implemented by performing gradient descent under a mixed canonical form.

\subsection{Gradient flow perspective}
The NGD perspective admits a gradient flow characterization. Let \(\mathcal{F} \colon \C^{2^n} \to \R\) be a loss function so that \(\mathcal{L}(\theta):= \mathcal{F}(q_{\theta})\) is the induced loss function on the parameter space.  
The following proposition links the natural gradient algorithm in \Cref{eqn: NGD for TN} to a discretization of a projected gradient flow under \(\mathcal{F}\).
\begin{prop}\label{prop: projected gradient flow}
    For any site \(i \in \{ 1, \ldots, n\}\), we let \(S_{i} = \{(\delta G_{k})_{k = 1}^{n} \mid \delta G_{k} = 0 \, \forall \, k \not = i\}\) and we consider
    \[
    \delta \theta_t \gets \argmin_{\delta\theta \in S_{i}} \left<\nabla_{\theta}\mathcal{L}|_{\theta = \theta_{t}}, \delta\theta \right> + \frac{1}{2}\eta \lVert q_{\theta_{t} + \delta \theta} - q_{\theta_{t}} \rVert_{F}^2.
    \]
    In other words, we let \(\delta \theta_t\) be the solution to the minimization task in \Cref{eqn: NGD for TN} from only changing the tensor component at \(G_i\). 
    Then, one has
    \[
    q_{\theta_{t}+\delta\theta_t} = q_{\theta_{t}} - \eta^{-1}\Pi_{i}\left(\gradientat{\nabla_{q}\mathcal{F}}{q = q_{\theta_{t}}}\right),
    \]
    where \(\Pi_{i}\) denotes the projection onto the tangent space of varying \(q_{\theta}\) at the \(i\)-th tensor component. 
\end{prop}

\Cref{prop: projected gradient flow} is directly related to the single site update step in \Cref{alg:1}. In the setting of \Cref{prop: projected gradient flow}, one sees that the continuous limit of taking \(\eta \to \infty\) leads to the ODE 
\[
\dot{q} = -\Pi_{i}\left(\nabla_{q}\mathcal{F}\right),
\]
which is indeed a projected gradient flow based on the loss function \(\mathcal{F}\) on the tensor space. 

One can also derive \Cref{alg:1} under a gradient flow perspective. To approximately simulate the gradient flow \(\dot{q} = -\nabla_{q}\mathcal{F}\) on an MPS ansatz, one can consider an ODE by the following equation
\begin{equation}\label{eqn: projected GF cont}
    \dot{q} = -\sum_{i = 1}^{n}\Pi_{i}\left(\nabla_{q}\mathcal{F}\right).
\end{equation}
One can see that the site update schedule \Cref{alg:1} is exactly derived by using an operator splitting on \Cref{eqn: projected GF cont}. The procedure in \Cref{alg:1} is then a forward Euler method. We remark that this perspective is similar to the time-dependent variational principle (TDVP) \cite{haegeman2011time, haegeman2016unifying} in MPS literature.

\subsection{Numerical experiment with NGD}\label{sec: NGD numerics}

We demonstrate that the proposed NGD method allows the training dynamics to avoid the local minima issue in MPS tomography. 

\begin{figure}
    \centering
    \includegraphics[width=0.9\columnwidth]{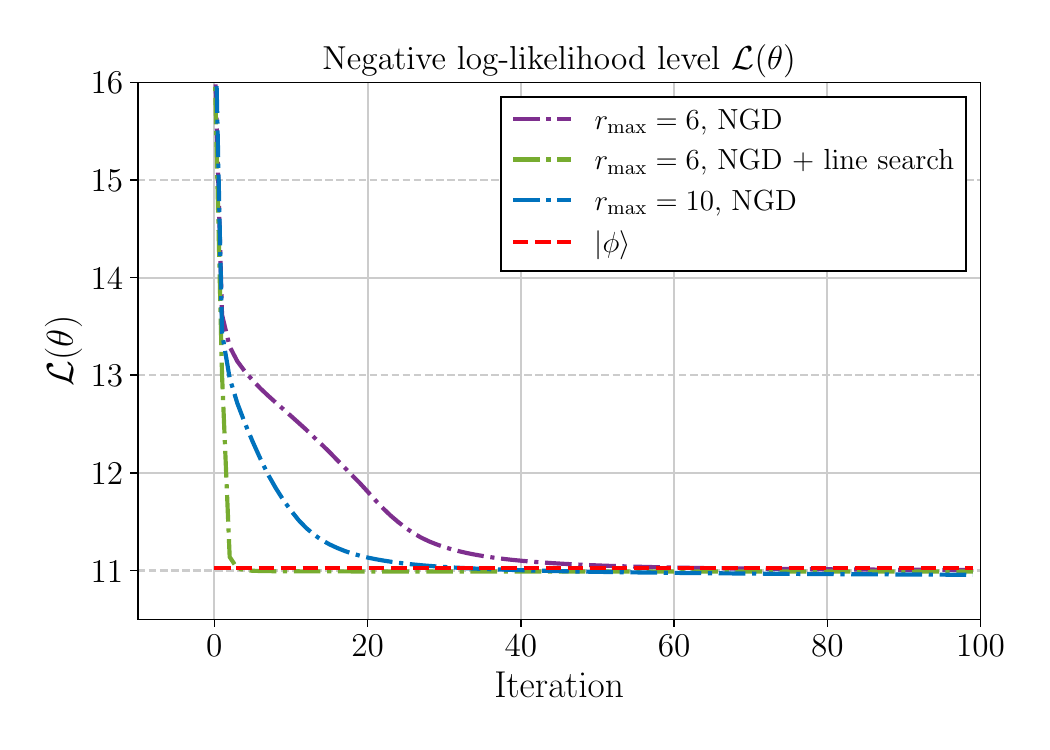}
    \caption{Performance of NGD for MPS tomography algorithm on the ground state of the periodic TFIM model in \Cref{eqn: TFIM}. The models are initialized randomly. The case where \(r_{\mathrm{max}}= 10\) with line search also converges rapidly and is omitted for simplicity.}
    \label{fig:result of NGD with MPS tomography}
\end{figure}

For the first example, we take the problem setting of \Cref{sec: experiment setting TFIM}. The result for applying the NGD method is illustrated in \Cref{fig:result of NGD with MPS tomography}. We see that all of the trained MPS models reach the NLL level of the ground state. Moreover, one can see that line search allows the MPS model to converge in only a few iterations.

\begin{figure}
    \centering
    \includegraphics[width=0.9\columnwidth]{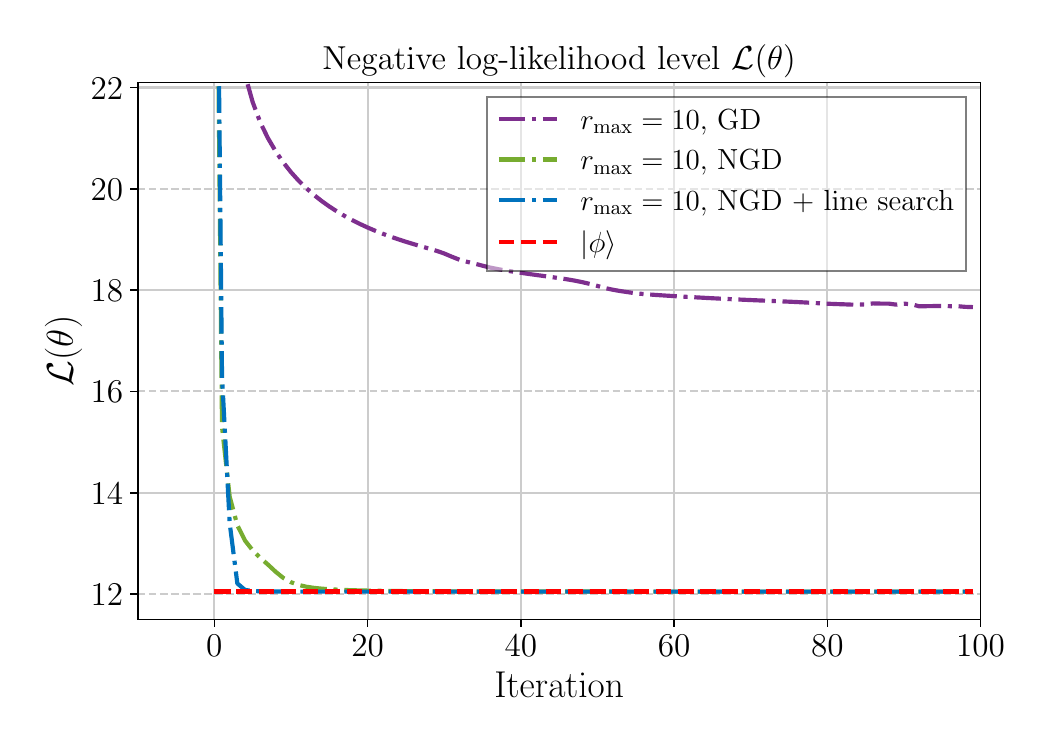}
    \caption{Performance of NGD for MPS tomography algorithm on the ground state of the periodic Heisenberg model in \Cref{eqn: heisenberg}. The models are initialized randomly. We choose the MPS model of maximal bond dimension \(r_{\mathrm{max}} = 10\), which is chosen according to the sample size to prevent overfitting.}
    \label{fig:result of NGD with MPS tomography for Heisenberg}
\end{figure}

For the second example, we consider a QST task for the ground state of the 1D antiferromagnetic Heisenberg model. The Hamiltonian of the model is
\begin{equation}\label{eqn: heisenberg}
    H = \sum_{(i,j) \in \mathrm{cycle}(n)} \left(\sigma_{i}^{X} \sigma_{j}^{X}+\sigma_{i}^{Y} \sigma_{j}^{Y}+\sigma_{i}^{Z} \sigma_{j}^{Z}\right),
\end{equation}
and the ground state is obtained by running DMRG with a maximal internal bond dimension \(r_{\mathrm{max}} = 40\).
We take \(n = 20\) for the system size, and we perform the maximum likelihood training based on \(B = 60000\) random Pauli measurements.  The result for applying the NGD method is illustrated in \Cref{fig:result of NGD with MPS tomography for Heisenberg}. The NGD method can quickly converge to the optimal NLL level, and NGD with line search can converge in a few iterations. In particular, the experiment shows that the training inefficiency of GD also occurs for antiferromagnetic models.

\begin{figure*}
    \centering
    \includegraphics[width=1.8\columnwidth]{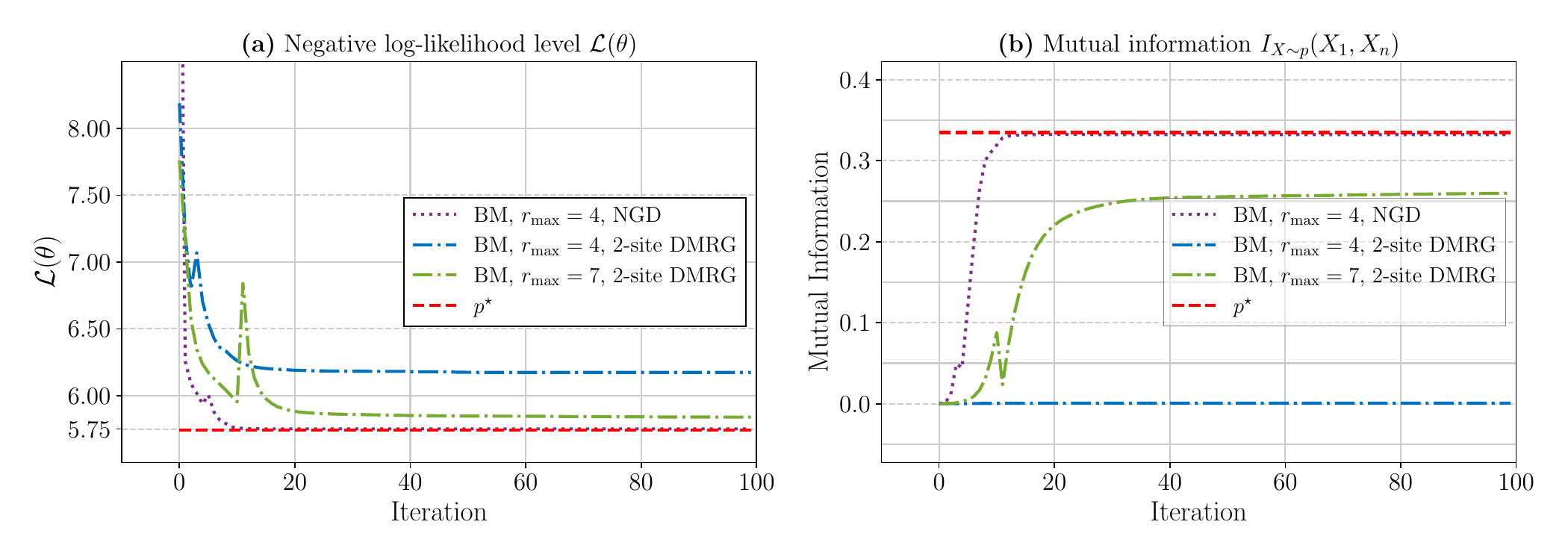}
    \caption{Performance of NGD for the Born machine algorithm on the periodic spin system model in \Cref{eqn: periodic model}. The models are initialized randomly. The 2-site DMRG method refers to the algorithm used in \cite{han2018unsupervised}, and is discussed in \Cref{appendix: NGD versus DMRG}. }
    \label{fig:result of NGD with BM}
\end{figure*}

For the third example, we show that NGD can address the causality trap in the Born machine. We apply NGD to the experiments of \Cref{sec: experiment setting}. The result is illustrated in \Cref{fig:result of NGD with BM}. We also compare the NGD method with the result of the training algorithm introduced in \cite{han2018unsupervised}. We refer to the algorithm in \cite{han2018unsupervised} as the 2-site DMRG method, and we defer a detailed discussion on this algorithm to \Cref{appendix: NGD versus DMRG}. One can see that the local minima issue occurs for the 2-site DMRG method. Moreover, the 2-site DMRG model also exhibits the causality trap when \(r_{\mathrm{max}} = 4\). In contrast, the NGD method can converge to the optimal NLL level.

\begin{figure*}
    \centering
         \includegraphics[width=1.8\columnwidth]{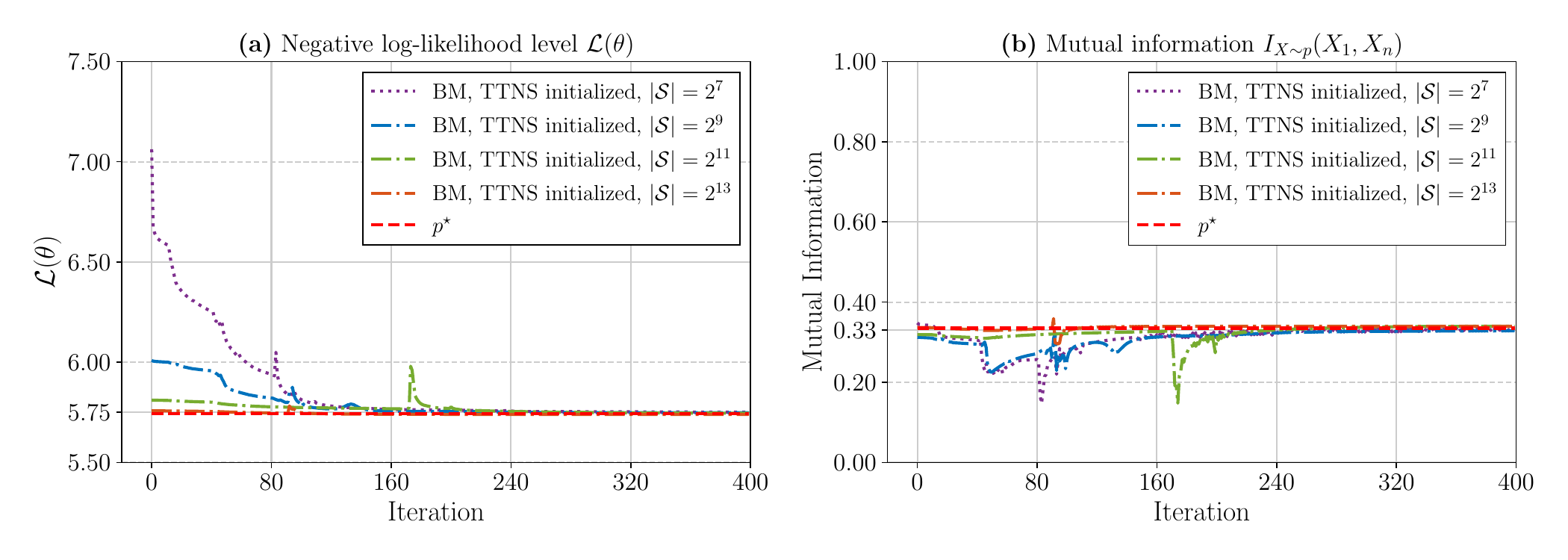}
    \caption{Performance of Born machine algorithm for the periodic spin system model. The models are initialized with TTNS-Sketch models from \Cref{alg: TT-cross fitting}. \(|\mathcal{S}|\) is the number of samples received by the TTNS-Sketch model. A warm initialization receiving only \(|\mathcal{S}| = 128\) samples is sufficient for the gradient descent algorithm to reach the optimal NLL level.}
    \label{fig:result of GD with TTNS}
\end{figure*}

\section{Avoiding causality trap with TTNS-Sketch initialization}\label{sec: TTNS initialization}

As is common in nonconvex optimization, one can avoid local minimum issues by a warm initialization that is close to the global minimum.
Our proposed strategy relies on a \emph{direct MPS ansatz}, which models a probability density by the following equation:
\begin{equation}\label{eq: TTNS-Sketch equation}
p_{\iota}(x) = \frac{q_{\iota}(x)}{W_{\iota}},
\end{equation}
where \(q_{\iota}\) is an MPS with tensor component \(\iota = (G_{k})_{k=1}^{n}\), and \(W_{\iota} = \sum_{z \in \{-1,1\}^n}q_\iota(z)\) is the normalization constant.

In contrast to the BM ansatz in \Cref{eq: BM equation}, the unique advantage of modeling distribution by a direct MPS ansatz is that it has a density estimation algorithm with a sample convergence guarantee. In particular, one can use the TTNS-Sketch algorithm \cite{tang2023generative}. The input to TTNS-Sketch is the samples \(\mathcal{T}\) drawn from \(p^{\star}\), and the output is an approximation of \(p^{\star}\) in a direct MPS ansatz. The TTNS-Sketch algorithm enjoys the following convergence guarantee (proof is in Appendix \ref{appendix: BM is MPS}):
\begin{prop}\label{prop: informal sample complexity TTNS-Sketch}
Let \(p^{\star}\) be an \(n\)-dimensional discrete distribution representable by a Born machine in \Cref{eq: BM equation}. Let \(\PTS\) denote the output of TTNS-Sketch algorithm after receiving \(B\) samples drawn from \(p^{\star}\). A sample size of \(B \ge \tilde{O}(\frac{n^2 + \epsilon  n}{\epsilon^2})\) ensures that \(\lVert p^{\star} - \PTS\rVert_{\infty} < \epsilon.\)
\end{prop}
Despite the guarantee, the direct MPS ansatz cannot supplant the BM ansatz as it can not guarantee to have only non-negative entries.
In this work, we propose to utilize the TTNS-Sketch output \(\PTS\) to form a warm initialization of a BM algorithm. Doing so allows one to leverage the convergence guarantee of a direct MPS ansatz and the positivity of a BM ansatz.

\subsection{Warm-start initialization of BM with TTNS-Sketch}
We explain the main idea of the proposed warm start procedure. After running the TTNS-Sketch algorithm in \cite{tang2023generative}, we obtain a direct MPS ansatz \(\PTS \approx p^{\star}\). Essentially, our strategy is to fit a BM ansatz against the square root of \(\PTS\). Utilizing the fact that accessing the entries of \(\PTS\) is efficient, one can obtain the BM ansatz by a simple MPS interpolation task. To accommodate different internal bond dimension specifications, we perform a postprocessing of the interpolation result, and the output of the fitting task is the warm BM initialization. The procedure is summarized in \Cref{alg: TT-cross fitting}.

We detail the steps taken in \Cref{alg: TT-cross fitting}. First, accounting for the possibly negative entry of \(\PTS\), we propose to use the TT-cross algorithm \cite{oseledets2010tt} to perform MPS interpolation with target function \(q_{>0}(z):= \sqrt{\max({0, \PTS(z)})}\). The output of the TT-cross is \(\theta_{\DMRG}\) so that \(p_{\theta_{\DMRG}} \approx \PTS\). 
Typically, the output of TT-cross is not of the specified internal bond dimension. Therefore, after obtaining the TT-cross output \(\theta_{\DMRG}\), we perform a MPS fitting task \(\theta_{\mathrm{init}} \gets \argmin_{\theta} \lVert q_{\theta} - q_{\theta_{\DMRG}}\rVert_{F}^2\), and \(\theta_{\mathrm{init}}\) is the warm initialization. In our case, the MPS fitting task is done by first performing a truncation of \(q_{\theta_{\DMRG}}\) to the specified internal bond dimension. Subsequently, we perform alternating least squares (ALS) to fit \(q_{\theta_{\DMRG}}\) until the task \(\min_{\theta} \lVert q_{\theta} - q_{\theta_{\DMRG}}\rVert_{F}^2\) reaches convergence.

\begin{algorithm}[H]
\caption{Warm start with TT-cross interpolation}
\label{alg: TT-cross fitting}
\begin{algorithmic}[1]
\REQUIRE TTNS-Sketch output \(\PTS\).
\REQUIRE TT-cross subroutine for MPS interpolation.

\STATE 
    \(\theta_{\DMRG} \gets \text{TT-cross}\left(\left(\sqrt{\PTS}\right)_{+}\right) \) \hfill Use TT-cross to fit square root of \(\PTS\). 
\STATE
    \(\theta_{\mathrm{init}} \gets \argmin_{\theta} \lVert q_{\theta} - q_{\theta_{\DMRG}}\rVert_{F}^2\). \hfill Post-processing of \(q_{\theta_{\DMRG}}\)

\RETURN $\theta_{\mathrm{init}}$ 
\end{algorithmic}
\end{algorithm}

\subsection{Numerical experiment with TTNS-Sketch initialization}

We demonstrate that the proposed warm initialization allows the training dynamics to avoid the causality trap. 
We take the problem setting of \Cref{sec: experiment setting}. The result for the performance of the warm initialization from Algorithm \ref{alg: TT-cross fitting} is illustrated in \Cref{fig:result of GD with TTNS}.
To assess the BM training under different qualities of warm initialization, we only use a subset \(\mathcal{S}\) of the total sample \(\mathcal{T}\) to obtain \(\PTS\), and we evaluate the training performance under different sample size \(|\mathcal{S}|\). \Cref{fig:result of GD with TTNS}(b) shows that all of the warm initialized models match \(p^{\star}\) in the mutual information on the variable pair \((X_1, X_n)\). Moreover, the high mutual information level in \Cref{fig:result of GD with TTNS}(b) suggests that the training dynamics of all of the models are quite far away from \(\pcausal\).

We draw two more conclusions from the warm initialization regarding the quality of the warm initialization. First, we see that cases of large sample size \(|\mathcal{S}|\) result in a high accuracy TTNS-Sketch model \(\PTS\), and the output of the warm initialization from Algorithm \ref{alg: TT-cross fitting} is already close to the optimal NLL level. Moreover, at \(|\mathcal{S}| = 2^{13}\), one can see that the initialized probability function \(p_{\theta_{\mathrm{init}}}\) is already sufficiently close to \(p^{\star}\).

Secondly, the accuracy requirement on \(\PTS\) for the warm initialization is quite mild. From \Cref{fig:result of GD with TTNS}, we see that the training is successful even if the warm initialization is from a TTNS-Sketch output \(\PTS\) obtained from only \(|\mathcal{S}| = 128\) samples. Our result suggests that a successful BM training does not require the warm initialization to be very close to the global minimum.

\section{Conclusion and outlook}\label{sec: conclusion}
This work studies trainability issues with MPS tomography when trained using standard gradient descent methods. We propose two effective solutions to avoid local minima based on a natural gradient algorithm and a warm start initialization. For the Born machine algorithm, we see that a high-quality warm initialization is already at the optimal NLL level. For practical examples of quantum state tomography, we see that the NGD method with line search allows rapid convergence to the optimal NLL level within just a few iterations. An open question is whether one can have a warm initialization subroutine for general MPS state tomography tasks for models with non-local interactions.

\bibliography{ref}

\begin{appendix}

\section{Characterization of the causality trap}\label{appendix: KL analysis}
We shall show that the causality trap can be characterized as an inability of a trained MPS model to capture boundary correlations.
We consider the NLL loss \(\NLL\) in the limit of sample size \(|\mathcal{T}| \rightarrow \infty\). Under this limit, the loss function reduces to the Kullback-Leibler (KL) divergence \(\KL{\cdot}{\cdot}\) as \[\NLL(\theta) \rightarrow \KL{p^{\star}}{p_\theta} + C,\] where \(C\) is a constant independent of \(\theta\). Define \(\theta^{\star}\) to be an MPS configuration so that \(p_{\theta^{\star}} = p^{\star}\). The NLL is minimized by taking \(\theta = \theta^{\star}\). Under this limit, the NLL gap is the KL divergence, as one has \[\mathcal{L}(\theta) - \mathcal{L}(\theta^{\star}) \rightarrow \KL{p^{\star}}{p_\theta}.\]

Let \(p(x)\) be a likelihood function and let \(x = (z,w)\) be a partition of the joint variable \(x\). We write \(p_{z}\) as the likelihood function for the marginal distribution of \(z\). We write \(p_{w \mid z = a}\) as the likelihood function for \(w\) condition on \(z = a\).
We decompose the joint variables \((x_{1}, \ldots, x_{n})\) into \((z,w)\), where \(z := (x_{1}, x_{n}), w := (x_{2}, \ldots, x_{n-1})\). 
The KL divergence between two generic distributions satisfies a chain rule, which leads to the following decomposition of \(\KL{p^{\star}}{p}\):
\begin{equation}\label{eqn: KL divergence decomposition}
\begin{aligned}
       \KL{p^{\star}}{p} = D_{z}(p) + D_{w \mid z}(p),
\end{aligned}
\end{equation}
where
\begin{align*}
&D_{z}(p) := \KL{p^{\star}_{z}}{p_{z}},
\\
&D_{w \mid z}(p) := \sum_{a}p^{\star}_{z}(a) \KL{p^{\star}_{w \mid z = a}}{p_{w \mid z = a}}.
\end{align*}

By direct computation, one sees that \(D_{w \mid z}(\pcausal) = 0\). Therefore, \(D_{z}(\pcausal)\) is the only term contributing to the NLL gap. 
Moreover, since \(X_{1}, X_{n}\) are close to being independent in \(\pcausal\), it follows that \(\pcausal_{(x_1, x_n)} \approx \pcausal_{x_{1}}\pcausal_{x_{n}} = p^{\star}_{x_{1}}p^{\star}_{x_{n}}\), where the second equality follows from the symmetry of \(p^{\star}\) and \(\pcausal\). Therefore, one has
\begin{equation}\label{eqn: causal model KL}
    \begin{aligned}
\KL{p^{\star}}{\pcausal}=&D_z(\pcausal) \\
= &\KL{p^{\star}_{(x_1, x_n)}}{\pcausal_{(x_1, x_n)}}\\
\approx &\KL{p^{\star}_{(x_1, x_n)}}{p^{\star}_{x_1}p^{\star}_{x_{n}}} \\
= & I_{X \sim p^{\star}}(X_{1}, X_{n}).
\end{aligned}
\end{equation}

The calculation in \Cref{eqn: causal model KL} shows that the NLL gap of \(\pcausal\) is essentially the mutual information for \((X_1, X_n)\) under \(p^{\star}\). Indeed, from \Cref{fig:result of GD}, one can see that the NLL gap for \(\pcausal\) is approximately \(0.33\), which coincides with the mutual information for \((X_1, X_n)\) in \(p^{\star}\). 

For practical cases with larger \(n\), it is no longer feasible to compare \(p_{\theta}\) with \(\pcausal\) with KL divergence or TV distance, as there is an \(O(2^n)\) cost in computing such metrics. The analysis we have given allows us to have a way to check the causality trap in practice. Formally, we characterize the causality trap as the setting in which a trained BM model \(p_{\theta}\) is a local minimum satisfying the following two conditions: (1) \(D_{w \mid z}(p_{\theta}) \approx D_{w \mid z}(\pcausal) = 0\), and (2) \(D_{z}(p_{\theta}) \approx D_{z}(\pcausal) \approx I_{X \sim p^{\star}}(X_{1}, X_{n})\). In other words, the causality trap occurs if the training algorithm succeeds in minimizing \(D_{w \mid z}(p_{\theta})\), but fails to minimize \(D_{z}(p_{\theta})\) beyond \(D_{z}(\pcausal)\).

One can verify the two given conditions of the causality trap by checking if \(p_{\theta}\) matches \(\pcausal\) in NLL level and mutual information. To see this, we illustrate how the plot in \Cref{fig:result of GD} implies that two stated conditions of the causality trap are met. \Cref{fig:result of GD}(b) shows that \((X_1, X_n)\) is approximately independent in any of the trained BM model \(p_{\theta}\), which shows that \(D_{z}(p_{\theta}) \approx D_{z}(\pcausal) \approx I_{X \sim p^{\star}}(X_{1}, X_{n})\). Then, \Cref{fig:result of GD}(a) shows that the NLL gap of \(p_{\theta}\) is approximately \(D_{z}(p_{\theta})\), which can only be true if \(D_{w \mid z}(p_{\theta}) \approx 0\).

\section{MPS training failure in a toy example}\label{sec: multi-linearity induced vanishing gradient}
We give a simple toy example to illustrate the potential training issues of the MPS ansatz under gradient descent.
Consider a family of MPS model \(q_{\theta}\) with tensor components \((c_{i}G_{i})_{i=1}^{n}\), where each \(G_{i}\) is fixed and each \(c_{i}\) is a scalar. 
In this case, the parameters are represented by \(\theta = (c_{1}, \ldots, c_{n})\). 
Let \(\mathcal{F} \colon \mathbb{C}^{2^n} \to \mathbb{R}\) be the loss function on \(q_{\theta}\). One can see that the tensor \(q_{\theta}\) only depends on \(x = \prod_{i = 1}^{n} c_{i}\), and therefore there exists a univariate loss function \(l\) so that \(\mathcal{F}(q_{\theta}) = l(\prod_{i = 1}^{n} c_{i})\). Without loss of generality, we assume that \(\lVert q_{\theta}\rVert_{F} = 1\) when \(x = \prod_{i = 1}^{n} c_{i} = 1\).

We suppose that \(l\) is strongly convex. In such a case, a reasonable optimization procedure is to perform gradient descent on \(x\) with learning rate \(\alpha\), where the update is given by
\[
x \gets x - \alpha\frac{dl}{dx}.
\]

First, we show that performing an NGD update step in \(c_{i}\) is equivalent to gradient descent in \(x\) with the same learning rate. Let \(\delta c_i\) be the update in \(c_{i}\). Let \(\delta q_{\theta}\) denote the associated update in \(q_{\theta}\). For \(\delta x =  \delta c_{k}\prod_{i \not = k}c_{i}\), one sees that \( \lVert \delta q_{\theta} \rVert = \lvert \delta x \rvert\). With a learning rate of \(\alpha\), the NGD step is done through the following formula:
\begin{equation}\label{eq: scalar TN dynamics}
    c_{i} \gets c_{i} + \argmin_{\delta c_{i}} \frac{\partial \mathcal{F}(q_{\theta})}{\partial c_{i}} \delta c_{i} + \frac{1}{2}\alpha^{-1} (\delta x)^2.
\end{equation}
One sees that \(\frac{\partial \mathcal{F}(q_{\theta})}{\partial c_{i}} \delta c_i = \frac{d l}{d x} \delta x\). Therefore, the resultant update to \(x\) follows the equation
\begin{equation*}
    x \gets x + \argmin_{\delta x} \frac{d l}{d x} \delta x + \frac{1}{2}\alpha^{-1} (\delta x)^2 = x - \alpha \frac{d l}{d x}.
\end{equation*}
Therefore, performing NGD in \(c_{i}\) with learning rate \(\alpha\) is equivalent to performing GD in \(x\) with the same learning rate.

In contrast, we show that performing a GD update step over \(c_i\) leads to instability. With a learning rate of \(\alpha\), the update to \(c_i\) is done by \(c_{i} \gets c_{i} - \alpha \frac{\partial \mathcal{F}(q_{\theta})}{\partial c_i},\)
and the resultant update to \(x\) is
\[
x \gets x - \alpha\left( \prod_{k \neq i}c_{k} \right)^2 \frac{dl}{dx}.
\]
One can see that performing GD in \(c_{i}\) with learning rate \(\alpha\) is equivalent to performing GD in \(x\) with learning rate \(\alpha \left(\prod_{i \not = k}c_{i}\right)^2\), which can be an exponentially large or exponentially diminishing learning rate for \(x\). Moreover, the formula shows that performing GD on a different site \(i\) leads to a different learning rate on \(x\).

Overall, the NGD method can better accommodate the multi-linear structure of the MPS ansatz. The given toy example illustrates the crucial observation that the NGD method tends to have fewer exploding or vanishing gradient problems, which allows for a more stable training performance.

\section{Proof of Proposition \ref{prop: projected gradient flow}}\label{appendix: proof of prop projected}

The proof is by direct computation. One has \(q_{\theta} \in \C^{2^n}\) and each tensor component \(G_{i}\) is a tensor of appropriate size and defined over \(\C\). In what follows, we split \(q_{\theta} \in \C^{2^n}\) into the real part and imaginary part, and we view \(q_{\theta}\) as a tensor in \(\R^{2*2^n}\). In the same way, we view each \(G_{i}\) as a tensor over \(\R\).

We view the tensor \(q_{\theta}\) and each tensor component \(G_{i}\) as having been flattened to a column vector of appropriate size.
For \(f \colon \R^{n} \to  \R^{m}\), let \(\gradientat{\frac{\delta f}{\delta x}}{x = a} \in \R^{m \times n}\) denote the Jacobian of \(f\) at \(x= a\).
Similarly, if  \(x = (z,w)\) is a partition of variables, then \(\gradientat{\frac{\delta f}{\delta z}}{x = a}\) is the submatrix of \(\gradientat{\frac{\delta f}{\delta x}}{x = a}\) constrained to columns corresponding to \(z\). 
If the codomain of \(f\) is \(\R\), the gradient satisfies \(\nabla_{x}f = \left(\frac{\delta f}{\delta x}\right)^{\top} \in \R^{n \times 1}\).

As a consequence of the multi-linearity of the MPS ansatz, for any \(\delta \theta \in S_{i}\), one has \[\gradientat{\frac{\delta q_{\theta}}{\delta G_{i}}}{\theta = \theta_{t} + \delta \theta} =
\gradientat{\frac{\delta q_{\theta}}{\delta G_{i}}}{\theta = \theta_{t}}.\]
Write 
\(
M = \gradientat{\frac{\delta q_{\theta}}{\delta G_{i}}}{\theta = \theta_{t}}\), 
\(L = \gradientat{\frac{\delta \mathcal{L}}{\delta G_{i}}}{\theta = \theta_t}\),
\(F = \gradientat{\frac{\delta \mathcal{F}}{\delta q} }{q = q_\theta}
\). Let \(\delta G_{i}\) be the update of \(\delta \theta \in S_{i}\) in the \(i\)-th tensor component. One has 
\[
q_{\theta_t  + \delta \theta} = q_{\theta_t} + M \delta G_{i}.
\]
One can write down the NGD update explicitly as a quadratic optimization:
\[
\begin{aligned}
\delta G_{i}^{\star} &= \argmin_{\delta G_{i}} L\delta G_{i} + \frac{1}{2}\eta \lVert q_{\theta_t  + \delta \theta} - q_{\theta_t}\rVert_{F}^2\\
&= \argmin_{\delta G_{i}} L\delta G_{i} + \frac{1}{2}\eta \delta G_{i}^{\top}M^{\top}M\delta G_{i}\\
&= \eta^{-1}\left(M^{\top}M\right)L^{\top}\\
&= \eta^{-1}\left(M^{\top}M\right)M^{\top}F^{\top},
\end{aligned}
\]
where the last equality follows from the chain rule.
Let \(\delta \theta \in S_{i}\) be the update to the tensor component so that the \(i\)-th component update is \(\delta G_{i}^{\star}\). Then
\begin{equation}\label{eqn: GF WTS QED}
    q_{\theta_t  + \delta \theta} - q_{\theta_t} =M\delta G_{i}^{\star} = \eta^{-1}M\left(M^{\top}M\right)^{-1}M^{\top}F^{\top}.
\end{equation}
Note that \(\Pi_{i}:= M\left(M^{\top}M\right)^{-1}M^{\top} = MM^{\dagger}\) is the projection onto the span of \(M\), and \(F^{\top} = \gradientat{\nabla_{q}\mathcal{F}}{q = q_{\theta_t}}\). Therefore, replacing \Cref{eqn: GF WTS QED} proves \Cref{prop: projected gradient flow} as is desired.

\section{Connection between NGD and DMRG-type algorithms}\label{appendix: NGD versus DMRG}

We shall show that the NGD update step in \Cref{alg:1} is equivalent to \Cref{alg:1-site dmrg}, which performs the GD update step in a mixed canonical form. As a result of the equivalence, one way to implement the NGD method is to apply gauge transformations. 

We prove that the update step in \Cref{alg:1-site dmrg} is equivalent to the NGD update step. Let \(\theta_{t}= (G_{k})_{k=1}^{n}\) be in a mixed canonical form centered at site \(i\). The update in \Cref{eqn:1-site dmrg} can be written by quadratic optimization:
\begin{equation}
    G_{i} = G_{i} + \argmin_{\delta G_{i}} \left<\nabla_{G_{i}}\mathcal{L}|_{\theta = \theta_t}, \delta G_{i}\right> + \frac{1}{2}\eta \lVert \delta G_{i} \rVert_{F}^2.
\end{equation}
Let \(\delta \theta \in S_{i}\) and let \(\delta G_{i}\) be its \(i\)-th tensor component. Because \(\theta_t\) is in the mixed canonical form centered at \(k\), it follows
\[\frac{1}{2}\eta \lVert \delta G_{k} \rVert_{F}^2 = \frac{1}{2}\eta \lVert q_{\theta_t} - q_{\theta_t + \delta \theta} \rVert_{F}^2.\]
Thus, the update to \(\theta_{t}\) by \Cref{eqn:1-site dmrg} is equivalent to \(\theta_t \gets \theta_t + \delta \theta\), where
\[
\delta \theta = \argmin_{\delta\theta \in S_{k}} \left<\nabla_{\theta}\mathcal{L}|_{\theta = \theta_{t}}, \delta\theta \right> + \frac{1}{2}\eta \lVert q_{\theta_{t} + \delta \theta} - q_{\theta_{t}} \rVert_{F}^2,
\]
which is the NGD update in \Cref{eqn:1-site dmrg}.

\begin{algorithm}[H]
\caption{1-site DMRG method update}
\label{alg:1-site dmrg}
\begin{algorithmic}[1]
\REQUIRE Loss function \(\mathcal{L}\).
\REQUIRE Current tensor component \(\theta_{t}\)
\REQUIRE Parameter \(\eta\)

\FOR{\(i = 1,\ldots, n\)}
\STATE 
Apply gauge transformation to \(\theta_{t}\) to a mixed canonical form centered at \(i\).
\STATE
\begin{equation}\label{eqn:1-site dmrg}
    G_{k} \gets G_{k} - \eta^{-1}\nabla_{G_{k}}\mathcal{L}|_{\theta = \theta_t}.
\end{equation}
\ENDFOR
\STATE
\(\theta_{t+1} \gets \theta_{t}\)
\RETURN $\theta_{t+1}$ 
\end{algorithmic}
\end{algorithm}

Moreover, the equivalence between NGD and \Cref{alg:1-site dmrg} facilitates a comparison between our NGD proposal with the training algorithm in \cite{han2018unsupervised} for BM. In \cite{han2018unsupervised}, performing sequential tensor component update with mixed canonical form is referred to as an algorithm inspired by the density matrix renormalization group (DMRG) algorithm. While our NGD method performs single tensor component updates, the algorithm in \cite{han2018unsupervised} performs tensor component updates by merging and splitting neighboring tensor components. Therefore, to simplify the discussion, we refer to \Cref{alg:1-site dmrg} as the 1-site DMRG method, and we refer to the training algorithm in \cite{han2018unsupervised} as the 2-site DMRG method.

\Cref{sec: NGD numerics} shows that the 2-site DMRG method encounters the local minima issue with the \(r_{\mathrm{max}}\) case entering the causality trap, whereas the 1-site DMRG method successfully reaches the optimal NLL level. Therefore, \Cref{sec: NGD numerics} shows that 1-site DMRG has superior performance than 2-site DMRG in avoiding local minima issues for MPS tomography problems.

While this work does not focus on why 2-site DMRG encounters the local minima issue, we shall discuss the procedure of 2-site DMRG and discuss the plausible mechanism for the local minima issue during training. Let \(\theta_{t}= (G_{k})_{k=1}^{n}\) be the current tensor component and let \((i, i + 1)\) be a pair of neighboring sites. To update the tensor components in \((i, i+1)\), the first step in 2-site DMRG is the \emph{merging step}. In particular, one constructs an MPS with tensor component \(\tilde{\theta}_{t}= (G_{k})_{k=1}^{i-1} \cup (G_{i,i+1}) \cup (G_{k})_{k=i+2}^{n}\). The tensor component \(G_{i, i+1}\) is obtained by merging tensor components \(G_{i}\) and \(G_{i+1}\). In the general case where \(1 < i < i + 1 < n\), one writes
\begin{equation}\label{eqn: Merge step}
    \begin{aligned}
    &G_{i, i + 1}(\alpha_{i-1}, (x_i, x_{i+1}), \alpha_{i+1}) \\=&\sum_{\alpha_{i}}G_{i}(\alpha_{i-1}, x_{i}, \alpha_{i})G_{i+1}(\alpha_{i}, x_{i+1}, \alpha_{i+1}).
\end{aligned}
\end{equation}
Similarly, the cases where \(i = 0\) and \(i +1 = n\) follows likewise by respectively omitting the \(\alpha_{i}\) and \(\alpha_{i+1}\) index in \Cref{eqn: Merge step}.
After the merge, the parameter \(\tilde{\theta}_{t}\) still represents an MPS \(q_{\tilde{\theta}_{t}}\), and in particular one has \(q_{\tilde{\theta}_{t}} = q_{\theta_{t}}\). One has an induced loss function \(\tilde{\mathcal{L}}\) for which \(\tilde{\mathcal{L}}(\tilde{\theta}_{t}) = \mathcal{L}(\theta_t)\). 

The second step in 2-site DMRG is the \emph{optimization step}. Similar to \Cref{alg:1-site dmrg}, we apply gauge transformation to \(\tilde{\theta}_{t}\) to a mixed canonical form centered at \((i, i+1)\). Then, one performs the gradient descent by taking \[\tilde{G}_{i,i+1} = G_{i,i+1} - \eta^{-1}\nabla_{G_{i, i+1}}\mathcal{L}|_{\tilde{\theta} = \tilde{\theta}_t}.\]

The last step in 2-site DMRG is the \emph{truncation step}. To obtain an update to the tensor components \(G_i\) and \(G_{i+1}\), one performs a QR or SVD factorization to find the best rank \(r_{k}\) factorization of \(\tilde{G}_{i, i+1}\):
\begin{equation}\label{eqn: DMRG truncation}
\begin{aligned}
&\tilde{G}_{i, i+1}(\alpha_{i-1}, (x_{i}, x_{i+1}), \alpha_{i+1})\\ \approx  &\sum_{\alpha_{i} = 1}^{r_{i}}\tilde{G}_{i}(\alpha_{i-1}, x_{i}, \alpha_{i})\tilde{G}_{i+1}(\alpha_{i}, x_{i+1}, \alpha_{i+1}),
\end{aligned}
\end{equation}
and then the update to \(\theta_t\) is performed by letting \((G_{i}, G_{k+1}) \gets (\tilde{G}_{i}, \tilde{G}_{i+1})\). 

The 2-site DMRG method in \cite{han2018unsupervised} performs the aforementioned update steps for each neighboring pairs \((i, i+1)\) by iterating from \(i = 1\) to \(i = n - 1\). One likely explanation for the 2-site DMRG method encountering local minima issues is that the truncation step in the 2-site DMRG is not variational. The factorization \Cref{eqn: DMRG truncation} does not necessarily minimize the loss \(\mathcal{L}\), and it instead simply fits the tensor \(\tilde{G}_{i,i+1}\) in the sense of Frobenius norm. Therefore, it is possible for the update \((G_{i}, G_{i+1}) \gets (\tilde{G}_{i}, \tilde{G}_{i+1})\) to increase the loss \(\mathcal{L}\). Thus, one possible explanation for the local minima issue is that the update to \(\tilde{\theta}_t\) during the optimization step is offset by the subsequent truncation step. 

Finally, while NGD is equivalent to mixed canonical form optimization for MPS, we remark that the NGD interpretation generalizes to other tensor networks. For example, for positive MPS \cite{glasser2019expressive}, one cannot apply a gauge transformation to the tensor components, which will destroy the positivity structure of its tensor components. However, using the NGD interpretation, one can perform an optimization step without needing to perform gauge transformations. Similarly, the NGD method readily generalizes to other 1D tensor network ansatz such as LPS \cite{glasser2019expressive}. The use of NGD optimization in other tensor network structures is a promising future research direction.

\section{Proof of Proposition \ref{prop: informal sample complexity TTNS-Sketch}}\label{appendix: BM is MPS}

From \cite{tang2023generative}, it has been proven that the TTNS-Sketch algorithm can converge to a distribution \(p^{\star}\) with the rate in \Cref{prop: informal sample complexity TTNS-Sketch} as long as \(p^{\star}\) is representable by an MPS. Thus, \Cref{prop: informal sample complexity TTNS-Sketch} holds if a BM can be represented by an MPS. 
\Cref{lem: Direct MPS covers BM} below shows the representation hierarchy between the BM ansatz as in \Cref{eq: BM equation} and a direct MPS ansatz as in \Cref{eq: TTNS-Sketch equation}.

\begin{lemma}\label{lem: Direct MPS covers BM}(Proposition 2 of \cite{glasser2019expressive})
If a function \(p\) is representable by a Born machine or a locally purified state with maximal internal rank \(r\), then there exists a representation of \(p\) using an MPS with an internal bond dimension no larger than \(r^2\).
\end{lemma}

Moreover, \Cref{lem: Direct MPS covers BM} implies that the statement in \Cref{prop: informal sample complexity TTNS-Sketch} also holds if one instead assumes that \(p^{\star}\) is a locally purified state (LPS). We refer the interested reader to \cite{glasser2019expressive} for a detailed account of LPS. 

\end{appendix}

\end{document}